\documentclass{article}

\usepackage{arxiv_template/PRIMEarxiv}

\usepackage[utf8]{inputenc} 
\usepackage[T1]{fontenc}    
\usepackage{hyperref}       
\usepackage{url}            
\usepackage{booktabs}       
\usepackage{amsfonts}       
\usepackage{nicefrac}       
\usepackage{microtype}      
\usepackage{lipsum}
\usepackage{soul,xcolor}
\usepackage{blkarray}

\usepackage{caption}
\usepackage{multirow}
\usepackage{subcaption}
\usepackage{float}
\usepackage{amsmath,bm}
\usepackage{mathtools}
\usepackage{fancyhdr}       
\usepackage{graphicx}       
\graphicspath{{media/}}     
\usepackage{comment}
\title{On the Geometry Transferability of the Hybrid Iterative Numerical Solver for Differential Equations

}

\author{
  Adar Kahana, Enrui Zhang, Somdatta Goswami, George EM Karniadakis \\
  Division of Applied Mathematics \\
  Brown university \\
  Providence, RI 02912\\
  \texttt{\{adar\_kahana, enrui\_zhang, somdatta\_goswami, george\_karniadakis\}@brown.edu} \\
  \And
  Rishikesh Ranade, Jay Pathak \\
  CTO Office \\
  Ansys Inc \\
  Canonsburg, PA 15317 \\
  \texttt{\{rishikesh.ranade, jay.pathak\}@ansys.com} \\
}

\begin{document}
\maketitle

\begin{abstract}
The discovery of fast numerical solvers prompted a clear and rapid shift towards iterative techniques in many applications, especially in computational mechanics, due to the increased necessity for solving very large linear systems. Most numerical solvers are highly dependent on the problem geometry and discretization, facing issues when any of these properties change. The newly developed Hybrid Iterative Numerical Transferable Solver (HINTS) combines a standard solver with a neural operator to achieve better performance, focusing on a single geometry at a time. In this work, we explore the "T" in HINTS, i.e., the geometry transferability properties of HINTS. We first propose to directly employ HINTS built for a specific geometry to a different but related geometry without any adjustments. In addition, we propose the integration of an operator level transfer learning with HINTS to even further improve the convergence of HINTS on new geometries and discretizations. We conduct numerical experiments for a Darcy flow problem and a plane-strain elasticity problem. The results show that both the direct application of HINTS and the transfer learning enhanced HINTS are able to accurately solve these problems on different geometries. In addition, using transfer learning, HINTS is able to converge to machine zero even faster than the direct application of HINTS.
\end{abstract}


\section{Introduction}

Numerical simulations play a crucial role in scientific and engineering applications such as mechanics of materials and structures \cite{hughes2012finite,simo2006computational,hughes2005isogeometric,jing2002numerical,rappaz2003numerical, goswami2019adaptive,bharali2022robust}, bio-mechanics \cite{zhang2022g2,goswami2022neural}, fluid dynamics \cite{patera1984spectral,kim1987turbulence,cockburn2012discontinuous}, etc. The simulation approach is based on solving linear/nonlinear partial differential equations (PDEs). The efficiency and accuracy of a simulation approach is always comparable to conflict partners. This means that the quest for a more efficient numerical solver frequently results in lesser numerical accuracy. In engineering simulations, the main factor is to have an acceptable accuracy with feasible computational requirements for both memory utilization and computing time. Furthermore, the solution process must be stable and dependable. Therefore, determining an appropriate approach for the problem at hand is crucial, and usually determines the outcome of a simulation.

In traditional numerical solvers like the finite element method, we often reduce the complex differential equations defining the physical system to a system of linear equations of the form: $\left[\mathbf{K}\right]\{\boldsymbol{u}\} = \{\boldsymbol{f}\}$, where $\left[\mathbf{K}\right]$ is referred as the stiffness matrix; $\{\boldsymbol{f}\}$ is the force vector and $\{\boldsymbol{u}\}$ is the set of unknowns. A simple, yet not recommended way to solve for the set of unknown is to use the direct method by inverting the stiffness matrix and multiplying it with a force vector. However, the direct method fails in cases of large degrees of freedom (the stiffness matrix is in the order of a few millions) and/or the stiffness matrix is sparse. At this juncture, the iterative solvers come to the rescue. We start with an initial guess for $\{\boldsymbol{u}\}$, and gradually progress towards the true solution for $\{\boldsymbol{u}\}$. The solver iterates until a reasonable solution that meets the stopping criterion (typically an error tolerance value) is obtained. Iterative solvers are appropriate for large computing problems because they can frequently be parallelized more efficiently using algorithms. However, proper pre-conditioning of the stiffness matrix is a mandatory requirement. The solution's high oscillatory component can be solved efficiently using a dense mesh and few steps with a simple iterative method like Jacobi iteration or Gauss-Seidel (GS) method. It suffers from divergence for non-symmetric and indefinite systems over a coarse grid, as well as slow convergence associated with low-frequency eigen modes, restricting its application for large scale linear systems.

Recent advances in deep learning in addition to the developments in computational power have provided the means to employ neural networks as efficient approximators for PDEs. Their compositional character differs from the traditional additive form of trial functions in linear function spaces, where PDE solution approximations are built using Galerkin, collocation, or finite volume approaches. Their computational parametrization via statistical learning and large-scale optimization approaches makes them increasingly suitable for solving nonlinear and high-dimensional PDEs. However, the neural network often learns the low-frequency eigen modes, and tends to avoid the high-frequency modes. This phenomenon referred to as spectral bias is observed in  numerous applications of neural networks.  

In one of the recent works \cite{zhang2022hybrid}, we proposed an efficient approach to integrate synergistically the iterative solvers with deep neural operators to exploit the merits of both the solvers in turn and to overcome the challenges of either of them. The approach acronym as \emph{HINTS}, improves the convergence of the solution across the spectrum of eigen modes by leveraging the spectral bias of a deep neural operator. As observed in the seminal work of HINTS, the solution is flexible with regards to the computational domain and is transferable to different discretizations. In this work, we investigate the transferability properties, i.e., the ``\emph{T}'' in HINTS, with respect to domain adaptation. Specifically, the information from a model trained on a specific domain (\textit{source}), is employed to infer the solution on a different but closely related domain (\textit{target}).

Additionally, we integrate HINTS with the seminal work of operator level transfer learning \cite{goswami2022deep} to improve the convergence rate of HINTS on the target domain. In this scenario, we use a small number of labelled data to fine tune the target model using samples from the target domain. The model is initialized with the learnt parameters of the source model and is trained under a hybrid loss function, comprised of a regression loss and the Conditional Embedding Operator Discrepancy (CEOD) loss, used to measure the divergence between conditional distributions in a Reproducing Kernel Hilbert Space (RKHS). The target model is trained only for the deeper layers, acknowledging the widely accepted fact that the shallow layers are responsible for capturing the more general features.

\begin{figure}[!ht]
    \centering
    \includegraphics[width=10cm]{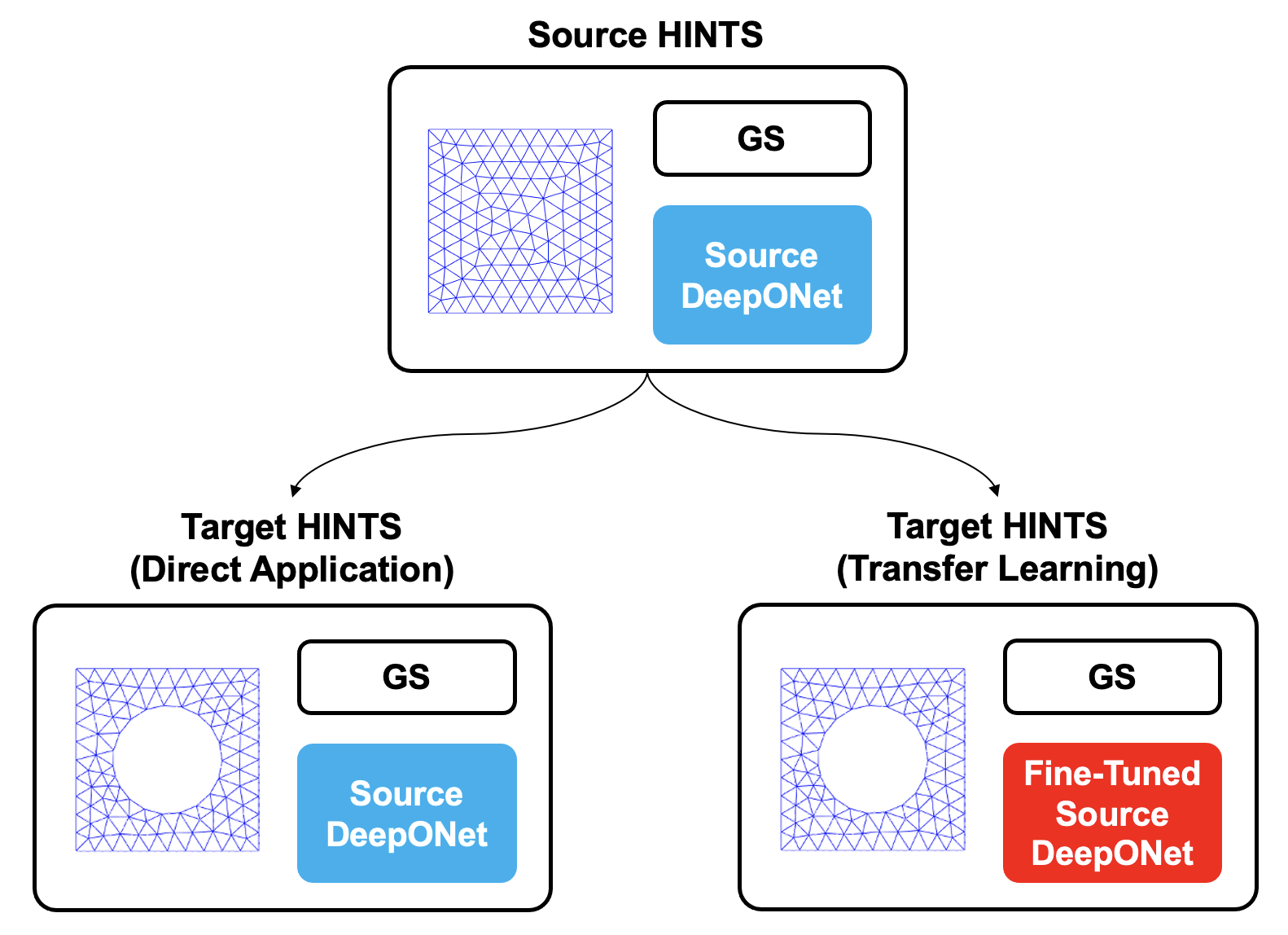}
    \caption{Diagram showing the proposed approaches. After building the source HINTS, one can use it directly on the target domain (left branch). In addition, we propose using transfer learning to fine tune the DeepONet of the source HINTS to get even better performance (right branch).}
    \label{fig:approaches}
\end{figure}

As a summary, we investigate the capability of HINTS trained on a source domain to operate on a target domain with the following two approaches: 
\begin{itemize}
    \item direct application of HINTS to an unseen target domain without retraining the DeepONet;
    \item usage of transfer learning to fine-tune the HINTS with limited data associated with the target domain.
\end{itemize}
This is illustrated in \autoref{fig:approaches}.

In this manuscript, we have considered the Darcy's model on a L-shaped domain (source), and the same domain with a circular or a triangular inclusion (target), and the linear elasticity model on a square domain (source) and a square domain with a circular inclusion (target). The remainder of the manuscript is organized as follows. In \autoref{sec:related_work}, we review the existing traditional numerical methods for solving linear systems, along with the recently popular deep learning solvers. This section also briefly covers the state-of-the-art neural operators, the Deep Operator Network (DeepONet), and a small discussion on the seminal work of HINTS and operator level transfer learning. In \autoref{approach}, we present the methodology of HINTS followed by the integration of the transfer learning approach. The experiments carried out to show the efficiency of domain adaptation with and without the transfer learning are presented in \autoref{sec:experiments}. Finally, we summarize our observations and provide concluding remarks in \autoref{sec:conclusion}.

\section{Related work}
\label{sec:related_work}

\textbf{Classical Numerical Solvers of ODEs/PDEs}. Over the last few decades, many studies have been conducted for developing advanced numerical solvers mainly for approximating the solutions of ODEs and PDEs. Classical methods are the Jacobi and GS methods, which were proposed in the 19th century. Since then, although many advanced algorithms have been proposed, most of the state-of-the-art solvers still use some versions of these two algorithms. The current leaders of the field of numerical solvers are the family of MultiGrid (MG) methods \cite{briggs2000multigrid,bramble2019multigrid}, where one uses a set of discretizations to solve the system. It is common to use Jacobi or GS smoothing inside the MG iterations. Therefore, constraints such as positive definiteness that apply for the Jacobi and GS algorithms, also affect MG as well. There are methods to overcome this, such as the shifted Laplacian method \cite{van2007spectral}, but these suffer from other disadvantages. The solvers compete for the lowest number of iterations needed for convergence, and also other properties such as physical time per iteration, parallelization capabilities, and more. However, in some cases the solvers may diverge (as an example, approximating the solution of an indefinite system), and solvers that are robust and can solve all types of problems are sought after.

\textbf{Machine Learning-Based Numerical Solvers}. Many authors have been investigating the use of AI when designing efficient numerical solvers. Some focus on creating an AI based solver for solving PDEs, replacing the numerical solvers. A notable example is the Physics-Informed Neural Networks \cite{raissi2019physics,goswami2020transfer}, where one trains a network to infer the solution of the PDE in a domain, without the need to assemble a linear system and solve, nor use a complex meshing algorithm. Another direction is to enhance numerical solvers using AI (which is the main focus of this paper). Most recent studies aim to replace components of the MG algorithms with AI based components, for example training a neural network to replace the restriction and prolongation operations \cite{moore2022learning,luz2020learning}. Others try to achieve a better performing preconditioner using learning \cite{gotz2018machine}.

\textbf{DeepONet}. Another notable advancement in the field of AI is the invention of operator learning methods \cite{goswami2022operator}. In contrast to standard Machine Learning (ML) tasks, where one seeks to approximate a function that can connect input data to output data, in operator learning one seeks a mapping between a family of functions and another family of functions that satisfies an operator, hence the name deep operator learning. Learning the operator enables many possibilities, fir example, after the network has been trained and the operator has been learned, one does not need to neither re-train nor solve the system again for new conditions or parameters, but rather infer the solution using the trained system for the new problem definition. This dramatically lowers the cost of solving PDE related problems. In addition, the learned operator does not depend on a discretization and can be used to infer the solution on any given discretization. Several operator learning methods have been proposed, including \cite{li2021fourier}. In this work, we focus on the DeepONet \cite{lu2021learning,goswami2022physics}, and use this for the operator learning.

\textbf{HINTS}. A new method to enhance numerical solvers using operator learning has been proposed under the name HINTS: Hybrid Numerical Iterative Transferable Solver \cite{zhang2022hybrid}. The idea of HINTS is to use an iterative solver such as Jacobi or GS, and replace some of the iterations with a DeepONet trained to receive the problem parameters and infer the solution. This DeepONet can also be used to receive the problem parameters and the residual at the current iteration and produce the correction term for the solution, so it can be applied in a similar way to the numerical solvers. Experiments show that using the existing numerical solvers, they tend to face slow convergence due to their difficulty to smooth the error for low frequency modes (where for the high modes they operate well), while the DeepONet excels in smoothing the low modes errors (but may fail for the high modes). Using both the numerical methods and the DeepONet, uniform convergence of all modes is achieved and the solvers converge to machine precision much faster. HINTS showed promising results on many tasks, and also a lot of potential for extensions, and in this paper we discuss a very important extension of HINTS that mechanics simulations may benefit from - the transferability of HINTS to new geometries and discretizations.

\textbf{Transfer Learning}. The idea of transfer learning is to leverage the parameters of a model trained with sufficiently large labelled dataset to infer information on a related task with few labelled datasets and minimal training. In \cite{goswami2022deep}, an operator level transfer learning was proposed to lower the computational costs of training a DeepONet (from Scratch) for related tasks. In the categorization of TL approaches, one popular classification is based on the consistency between the distributions of source and the target input (or feature) spaces and output (or label) spaces. The shift between the source and target data distributions is considered the major challenge in modern TL. One type of distribution shifts include conditional shift, where the marginal distribution of source and target input data remains the same while the conditional distributions of the output differ (i.e., $P(\mathbf{x}_s) = P(\mathbf{x}_t)$ and $P(\mathbf{y}_s\vert\mathbf{x}_s) \ne P(\mathbf{y}_t\vert\mathbf{x}_t)$) and covariate shift, where the opposite occurs (i.e., $P(\mathbf{x}_s) \ne P(\mathbf{x}_t)$ and $P(\mathbf{y}_s\vert\mathbf{x}_s) = P(\mathbf{y}_t\vert\mathbf{x}_t)$). In this work, we have employed the TL model proposed in \cite{goswami2022deep} to work on conditional distribution discrepancy between the domains. In this work, the authors have reported that the target model essentially requires a smaller dataset to fine tune, which can thus be done in significantly less time.

\section{Method} \label{approach}

Without any loss of generality, we consider a family of PDEs defined in a domain $\Omega$:
\begin{align}\label{eqn:diffeqn_general}
    \mathcal{L}_{\boldsymbol{x}}(\boldsymbol{u};k) &= \boldsymbol{f}(\boldsymbol{x}), \quad \boldsymbol{x} \in \Omega \\
    \mathcal{B}_{\boldsymbol{x}}(\boldsymbol{u}) &= \boldsymbol{g}(\boldsymbol{x}), \quad \boldsymbol{x} \in \partial\Omega \label{eqn:bound_general},
\end{align}
where $\mathcal{L}_{\boldsymbol{x}}$ is a differential operator, $\mathcal{B}_{\boldsymbol{x}}$ is a boundary operator, $k=k(\boldsymbol{x})$ parameterizes $\mathcal{L}_{\boldsymbol{x}}$,  $\boldsymbol{f}(\boldsymbol{x})$ and $\boldsymbol{g}(\boldsymbol{x})$ are the forcing terms, and $\boldsymbol{u} = \boldsymbol{u}(\boldsymbol{x})$ is the solution of the given PDE. With a well-trained DeepONet embedded, a HINTS is capable of solving for $\boldsymbol{u}$ corresponding to $k$ and $\boldsymbol{f}$, i.e., it captures the solution operator $\mathcal G$ of the family of PDEs specified by \autoref{eqn:diffeqn_general} and \autoref{eqn:bound_general}:
\begin{equation}
    \label{eqn:source_hints}
    \mathcal G: k, \boldsymbol{f} \mapsto \boldsymbol{u}\text{ s.t. \autoref{eqn:diffeqn_general} and \autoref{eqn:bound_general} are satisfied}.
\end{equation}
For detailed information on the implementation of HINTS, readers are suggested to refer to the work \cite{zhang2022hybrid}.

Herein, we first construct a source HINTS, i.e., a HINTS with DeepONet trained offline for \autoref{eqn:diffeqn_general} and \autoref{eqn:bound_general} defined in the source domain $\Omega=\Omega^S$. This is conducted by using $N_S$ labelled source data, $\mathcal D_S = \{k_i, \boldsymbol{f}_i, \boldsymbol{u}_i\}_{i = 1}^{N_S}$. In the DeepONet architecture, the branch network is a convolution neural network comprising of channels, each taking as input one of the mapping functions. The trunk network takes as inputs the spatial location of the points in the domain $\Omega^S$. This DeepONet is trained with a standard regression loss (relative mean squared error) to obtain the optimized weights and biases of the source network, $\boldsymbol{\theta}^S$. Next, we build HINTS by assembling the DeepONet and the numerical solver (e.g., GS). We discretize $\Omega^S$ with triangular elements. The simulation starts by assuming an initial solution of the dependent variable and at every iteration of HINTS it adopts either the pre-decided fixed relaxation method or the pre-trained DeepONet to approximate the solution. The iterative solver and DeepONet is chosen based on a pre-decided ratio. Among the available numerical iterative solvers, we have employed the GS approach.

With the source HINTS properly constructed, we now consider transferring it to become a target HINTS, i.e. a HINTS for solving the class of PDEs (\autoref{eqn:diffeqn_general} and \autoref{eqn:bound_general}) defined in a target domain $\Omega=\Omega^T$, with the following two methods:

\textbf{Direct Application: Use Source DeepONet Directly in Target HINTS}. To do this, the first approach is to directly apply the source HINTS for inference in the target domain $\Omega^T$. Specifically, when implementing the target HINTS (i.e., HINTS for the target domain), the trained DeepONet from the source HINTS is directly invoked in the workflow of the target HINTS.

\textbf{Transfer Learning: Fine-tune the DeepONet in Target HINTS}. In the second approach, we fine-tune the trained DeepONet from the source HINTS with a small number of labelled samples from the target domain. We generate $N_T$ labelled data, $\mathcal D_T = \{\tilde{k}_i, \tilde{\boldsymbol{f}}_i, \tilde{\boldsymbol{u}}_i\}_{i = 1}^{N_T}$ on the target domain, where $N_S \gg N_T$. In the examples presented in this work, $N_T \approx 0.01\times N_S$. We initialize a target model (having the same architecture as the source model), with the learnt parameters of the source model and fine tune the model by training the fully connected layers of the convolution neural network and the last layer of the trunk net under a hybrid loss function. The hybrid loss function reads as:
\begin{equation}
\label{eq:target-loss}
    \mathcal L(\theta^T) = \lambda_1 \mathcal L_r(\theta^T) + \lambda_2 \mathcal L_{\text{CEOD}}(\theta^T), 
\end{equation}
where $\lambda_1 = 1$ and $\lambda_2\gg\lambda_1$ are trainable coefficients, which determine the importance of the two loss components during the optimization process \cite{kontolati2022influence}. In \autoref{eq:target-loss}, while the first term is a standard regression loss, the second term ensures the agreement between the conditional distributions of the target data. For details on the construction of the $\mathcal L_{\text{CEOD}}(\theta^T)$, readers are suggested to refer to \cite{goswami2022deep}. While employing the HINTS algorithm for the target domain, the source DeepONet is replaced by the target DeepONet.

Even though the HINTS solution is transferable between related domains, and integrating the transfer learning approach would essentially mean an additional training time, we argue that the convergence rate of a transfer-learning integrated HINTS solution is in faster than its counterpart (see section \autoref{sec:conclusion} for statistical analysis). It is worth noting that the target model is trained with much less iterations and hence is very fast. We mention that there is potential for performing the transfer learning even better, and producing a more accurate target DeepONet, but in this work we did not aim for achieving this. We focus on showing that the transfer learning mechanism works and produces even faster convergence for HINTS.

\section{Computational experiments}
\label{sec:experiments}

In this section, we explore our method using two benchmark problems. The first problem involves the flow in heterogeneous porous media (Darcy's model) on a two-dimensional L-shaped domain (source). The target domains considered in this case are the same L-shaped domain with a circular and a triangular cutout. In the next problem, we have considered a thin rectangular plate (source) subjected to in-plane loading that is modeled as a two-dimensional problem of plane strain elasticity. The target domain considered in this case is the same rectangular plate with a central circular cutout. In both the examples, the DeepONet is trained using the Adam optimizer \cite{kingma2014adam}. The implementation has been carried out using the \texttt{PyTorch} framework \cite{paszke2019pytorch}. For both the examples we initialize the DeepONet parameters using Xavier initialization. Details on the data generation, network architecture, such as number of layers, number of neurons in each layer, and the activation functions are provided with each example.

\subsection{Darcy flow}
\label{subsec:darcy}

\begin{figure}[!ht]
     \centering
     \begin{subfigure}[t]{0.32\textwidth}
         \centering
         \includegraphics[width=1\textwidth]{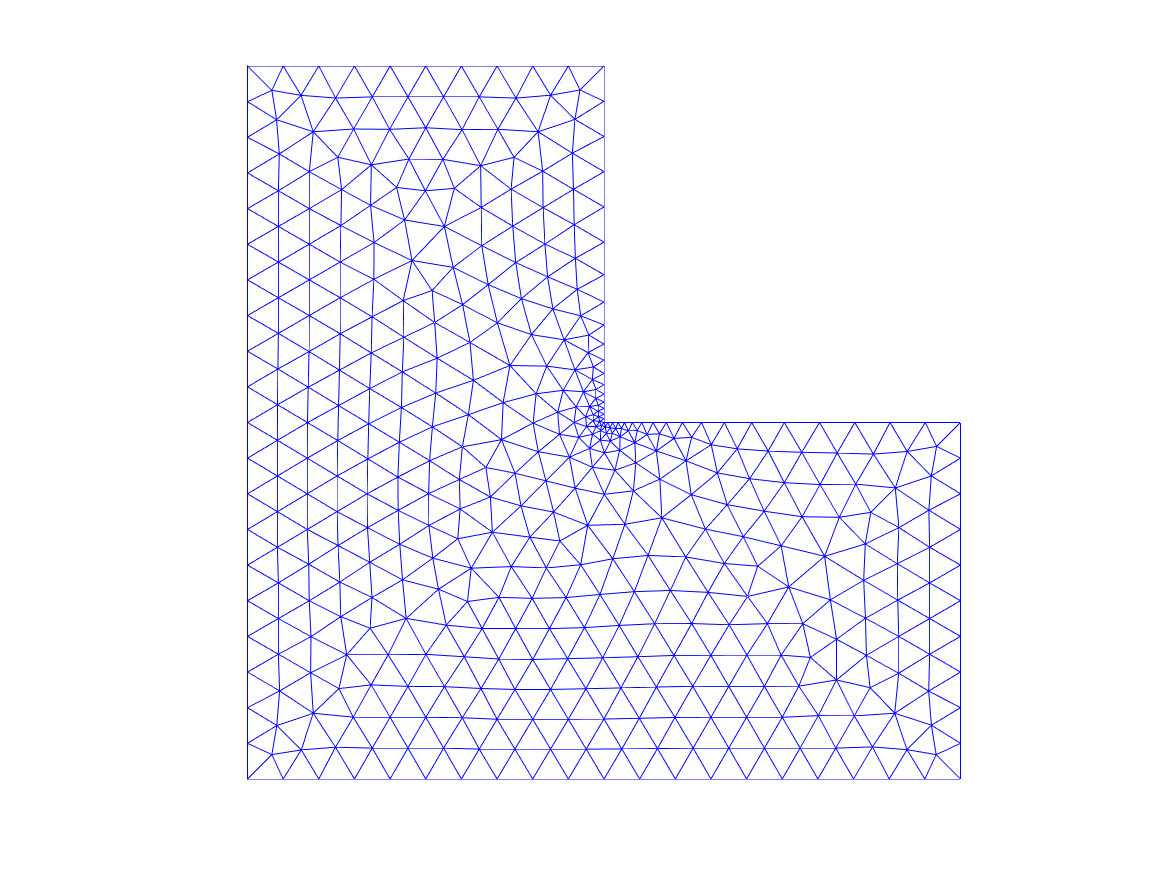}
         \caption{L-shaped domain.}
     \end{subfigure}
     \hfill
     \begin{subfigure}[t]{0.32\textwidth}
         \centering
         \includegraphics[width=1\textwidth]{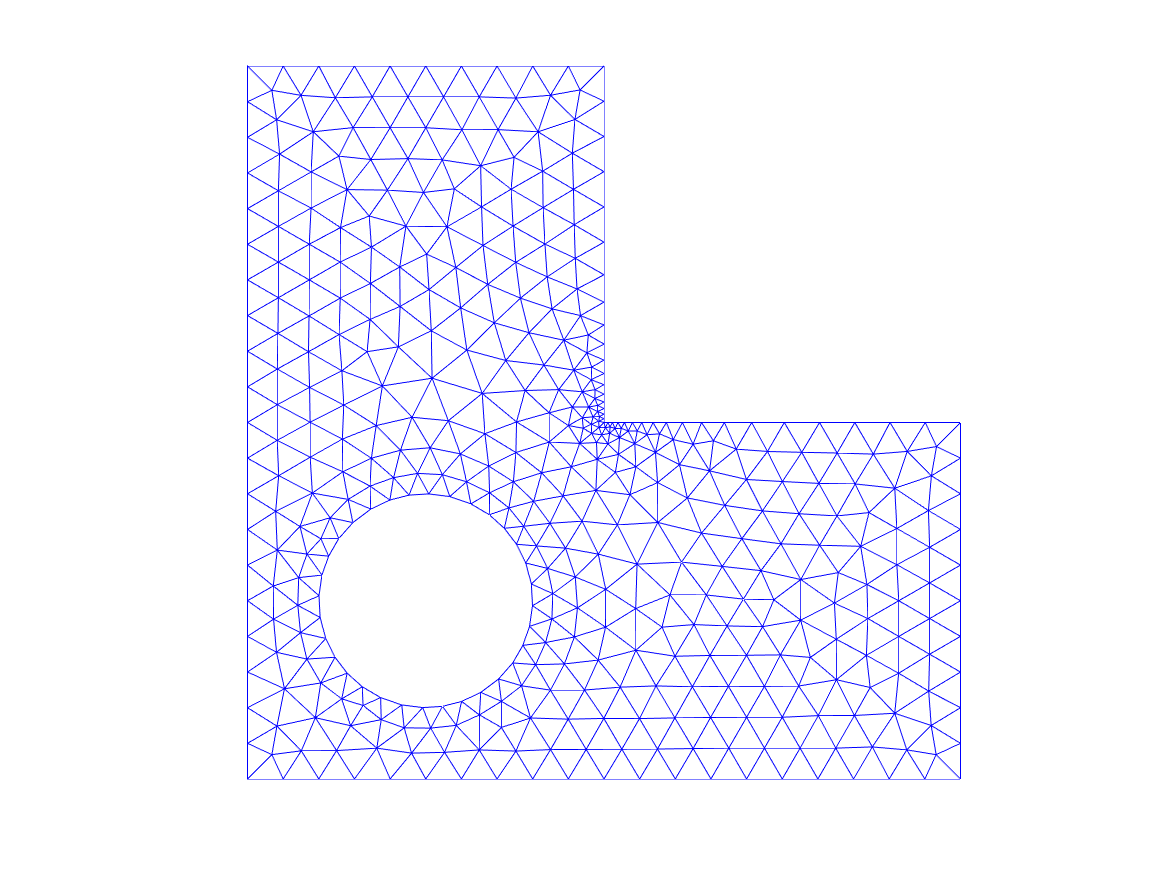}
         \caption{L-shaped domain with a circular hole.}
     \end{subfigure}
     \hfill
     \begin{subfigure}[t]{0.32\textwidth}
         \centering
         \includegraphics[width=1\textwidth]{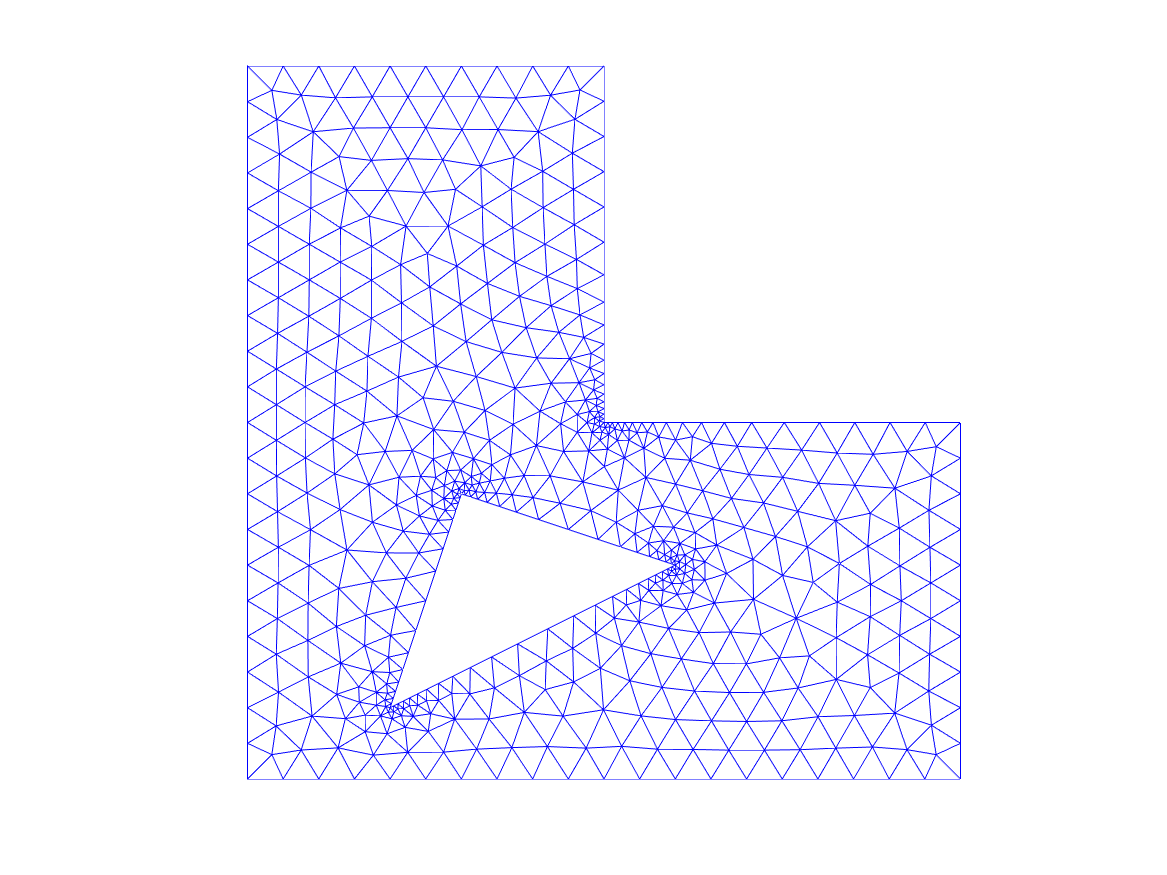}
         \caption{L-shaped domain with a triangular hole.}
     \end{subfigure}
     \hfill
     \caption{Mesh discretization for the geometries considered in the Darcy's problem: (a) source domain, (b-c) target domains considered for domain adaptation.}
     \label{fig:domains}
\end{figure}

In the first example we consider the Darcy's flow on a L-shaped domain. 
The problem is defined on a L-shaped domain (source) as:
\begin{align}
    \label{eqn:poisson_PDE}
    \nabla\cdot\left(k(\boldsymbol{x})\nabla u(\boldsymbol{x})\right) + f(\boldsymbol{x}) &= 0,\;\;\; \bm x = (x,y)\in\Omega^S_L:=(0,1)^2\backslash[0.5,1)^2 \\
    u(\boldsymbol x) &= 0, \quad \boldsymbol x\in\partial\Omega^S_L,
\end{align}
where $k(\boldsymbol{x})$ is a spatially varying hydraulic conductivity, $u(\boldsymbol x)$ is the hydraulic head, and $f(\boldsymbol x)$ is a spatially varying force vector. We use triangular elements to discretize the L-shaped domain, $\Omega^S_L$. The spatial discretization is shown in \autoref{fig:domains}(a).

\bigbreak
\noindent
\textbf{Training Source HINTS}

In this example, the goal is to learn the operator of the system in \autoref{eqn:poisson_PDE}, which maps the random conductivity field and the random force vector to the output hydraulic head, i.e., $\mathcal{G}_\theta: k(\boldsymbol x), f(\boldsymbol x) \rightarrow h(\boldsymbol x)$, where $\mathcal{G}_\theta$ is the solution operator. To generate multiple samples of the conductivity fields and force vectors for training the source DeepONet, we describe $k(\boldsymbol x)$ and $f(\boldsymbol{x})$ as stochastic processes, the realizations of which are generated using a Gaussian Random Field (GRF), with a standard deviation of $0.3$ and $0.1$ for $k(\boldsymbol x)$ and $f(\boldsymbol x)$, respectively and a correlation length of $0.1$ for both the processes. To train the source DeepONet, we generate $N_S = 51\small{,}000$ samples as the labeled dataset of random fields and the corresponding responses, and an additional $N_S^{test} = 9\small{,}000$ samples to test the model. 

The branch network of the DeepONet is a combination of a CNN (input dimension $31 \times 31$, number of channels $[2, 40, 60, 100, 180]$, kernel size $3\times 3$, stride $2$, the number of channels of the input 2 comes from the concatenation of $K(\boldsymbol x)$ and $f(\boldsymbol x))$ and a fully-connected network (dimensions $[180, 80, 80]$). The dimension of the trunk network (fully-connected network) is $[2, 80, 80, 80]$. We train the DeepONet for $25\small{,}000$ epochs with a fixed learning rate of $1e-4$, and a mini-batch size of $10\small{,}000$. For the branch network we employ the \texttt{ReLU} activation function and for the trunk network, we use \texttt{Tanh} activation. The last layers of both the branch and the trunk networks have linear activation functions. The source model converges with a mean relative error of $4\%$ on the test dataset.

\begin{figure}[htbp!]
     \centering
     \begin{subfigure}[t]{0.35\textwidth}
         \centering
         \includegraphics[width=\textwidth]{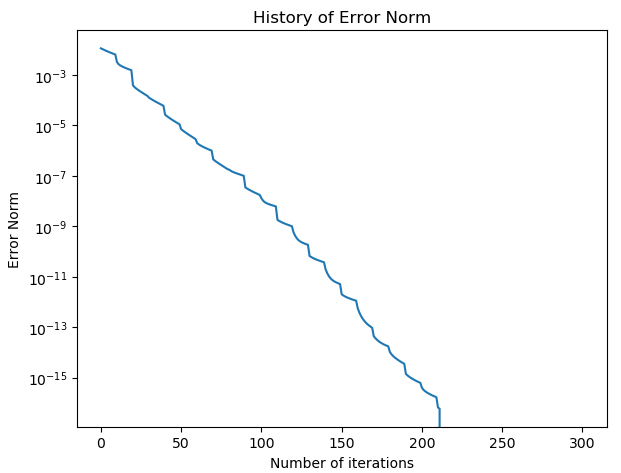}
         \caption{The error norm per iteration.}
         \label{fig:HINTS_Poisson_error_norm}
     \end{subfigure}
     \begin{subfigure}[t]{0.35\textwidth}
         \centering
         \includegraphics[width=\textwidth]{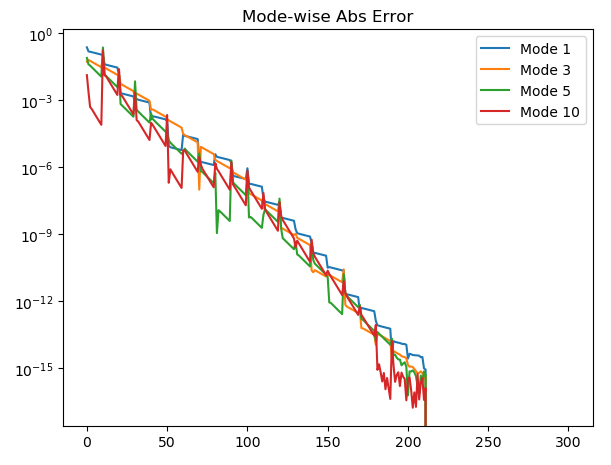}
         \caption{The mode errors.}
         \label{fig:HINTS_Poisson_modes}
     \end{subfigure}
     \begin{subfigure}[t]{0.35\textwidth}
         \centering
         \includegraphics[width=\textwidth]{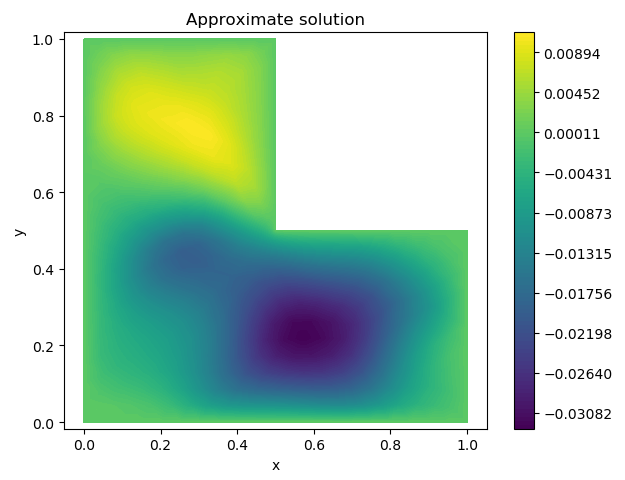}
         \caption{The approximate solution at the last step.}
         \label{fig:HINTS_Poisson_sol}
     \end{subfigure}
     \begin{subfigure}[t]{0.35\textwidth}
         \centering
         \includegraphics[width=\textwidth]{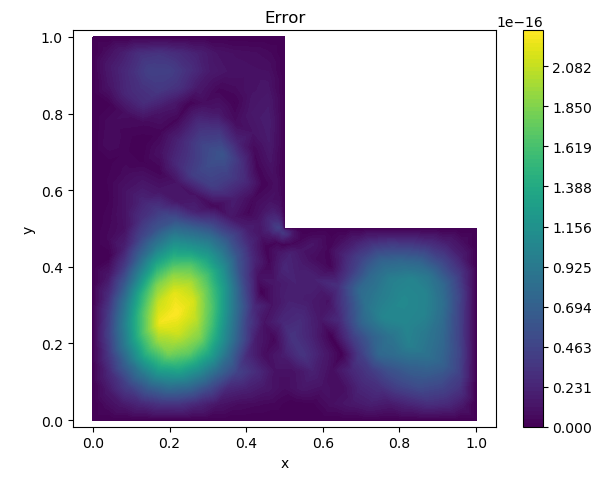}
         \caption{The error at the last iteration of HINTS.}
         \label{fig:HINTS_Poisson_error}
     \end{subfigure}
     \hfill
     \caption{Employing HINTS to solve the Darcy's problem on a L-shaped domain. It is interesting to note that in (d) the scale shows convergence on the error to machine zero at the last iteration of the proposed solver.}
     \label{fig:HINTS_Poisson}
\end{figure}

The trained DeepONet is plugged into the HINTS algorithm and the HINTS iterations are executed until the error of the solution reaches machine zero. An example of such iterative solution is given in \autoref{fig:HINTS_Poisson}. When inspecting the solution, we focus mainly on the convergence efficiency of the solution to machine zero. For that, we observe the norm of the error per iteration, as demonstrated in \autoref{fig:HINTS_Poisson_error_norm}. The founding idea of HINTS is based on the uniform convergence of all the eigen modes, which is shown in \autoref{fig:HINTS_Poisson_modes}. To illustrate the solution of the problem at hand, we show \autoref{fig:HINTS_Poisson_sol}, which is the solution at the last iteration (after convergence). We also show the error between the approximate solution and the exact solution (obtained by directly solving the system instead of using an iterative method) in \autoref{fig:HINTS_Poisson_error}. Note that the scale of the error is machine precision limit, which means the desired convergence has been achieved.

Next, to investigate the domain adaptation capabilities, we consider the following two tasks:
\begin{itemize}
    \item \textbf{Task1:} From a L-shaped domain to a L-shaped domain with a circular cutout. The target domain is defined as: $\Omega^{T1}_L = \Omega_L \backslash \{(x, y) | (x - 0.25) ^ 2 + (y - 0.25) ^ 2 \leq 0.15\}$.
    \item \textbf{Task2:} From a L-shaped domain to a L-shaped domain with a triangular cutout. The target domain is defined as: $\Omega^{T2}_L = \Omega_L \backslash \{(x, y) | (x, y) \in \bigtriangleup((0.2, 0.1), (0.6, 0.4), (0.3, 0.4))\}$. 
\end{itemize}
The discretization of the target domains are shown in \autoref{fig:domains} (b) and (c). While the circular cutout has a smooth boundary and is easier to approximate, the triangular cutout has locations of singularity, and hence imposes a more challenging scenario for domain adaptation. 

\begin{figure}[htbp!]
     \centering
     \begin{subfigure}[t]{0.49\textwidth}
         \centering
         \includegraphics[width=\textwidth]{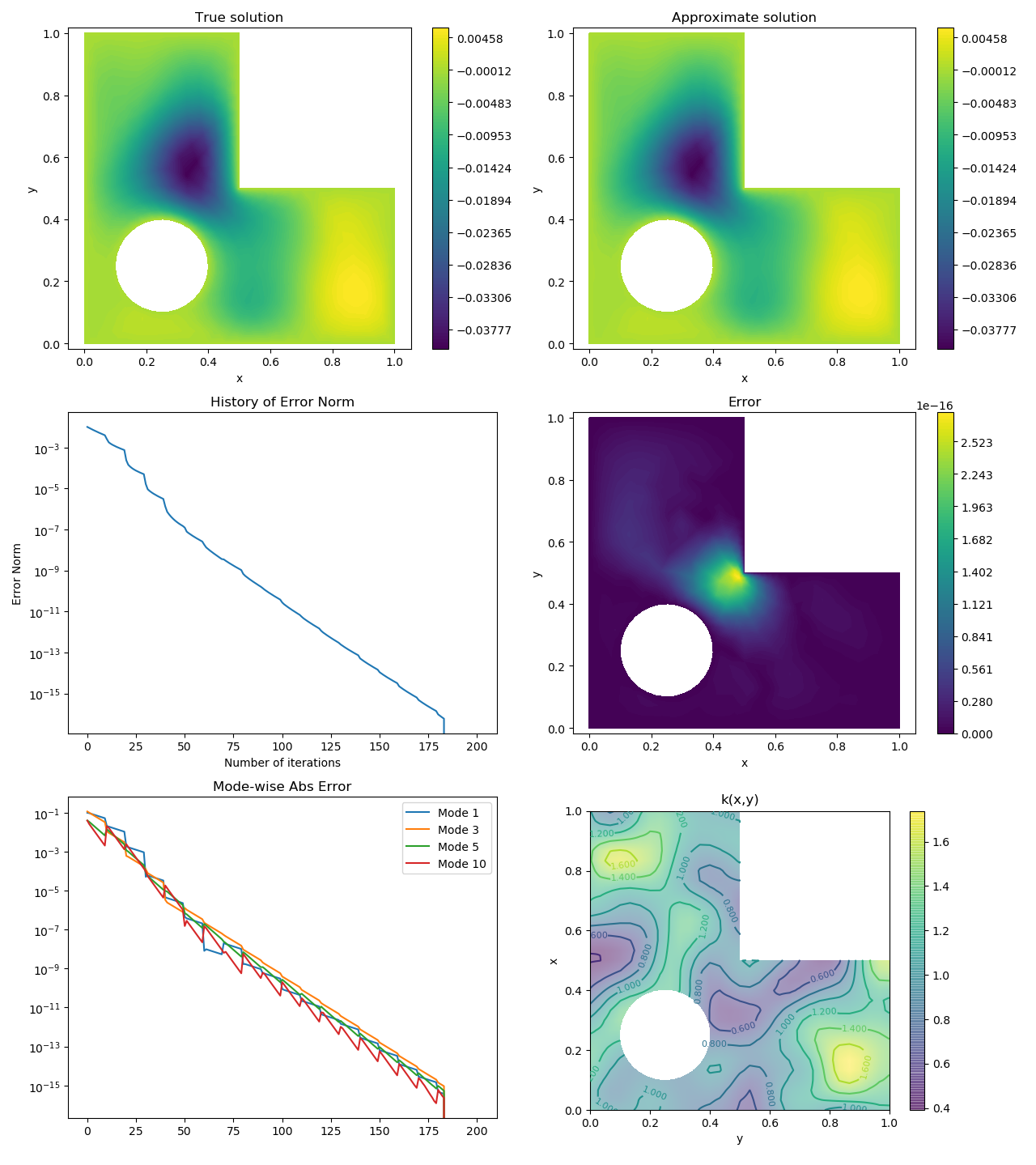}
         \caption{HINTS trained on $\Omega_L$, solving on $\Omega_L^{T1}$.}
         \label{fig:HINTS_circle}
     \end{subfigure}
     \hfill
     \begin{subfigure}[t]{0.49\textwidth}
         \centering
         \includegraphics[width=\textwidth]{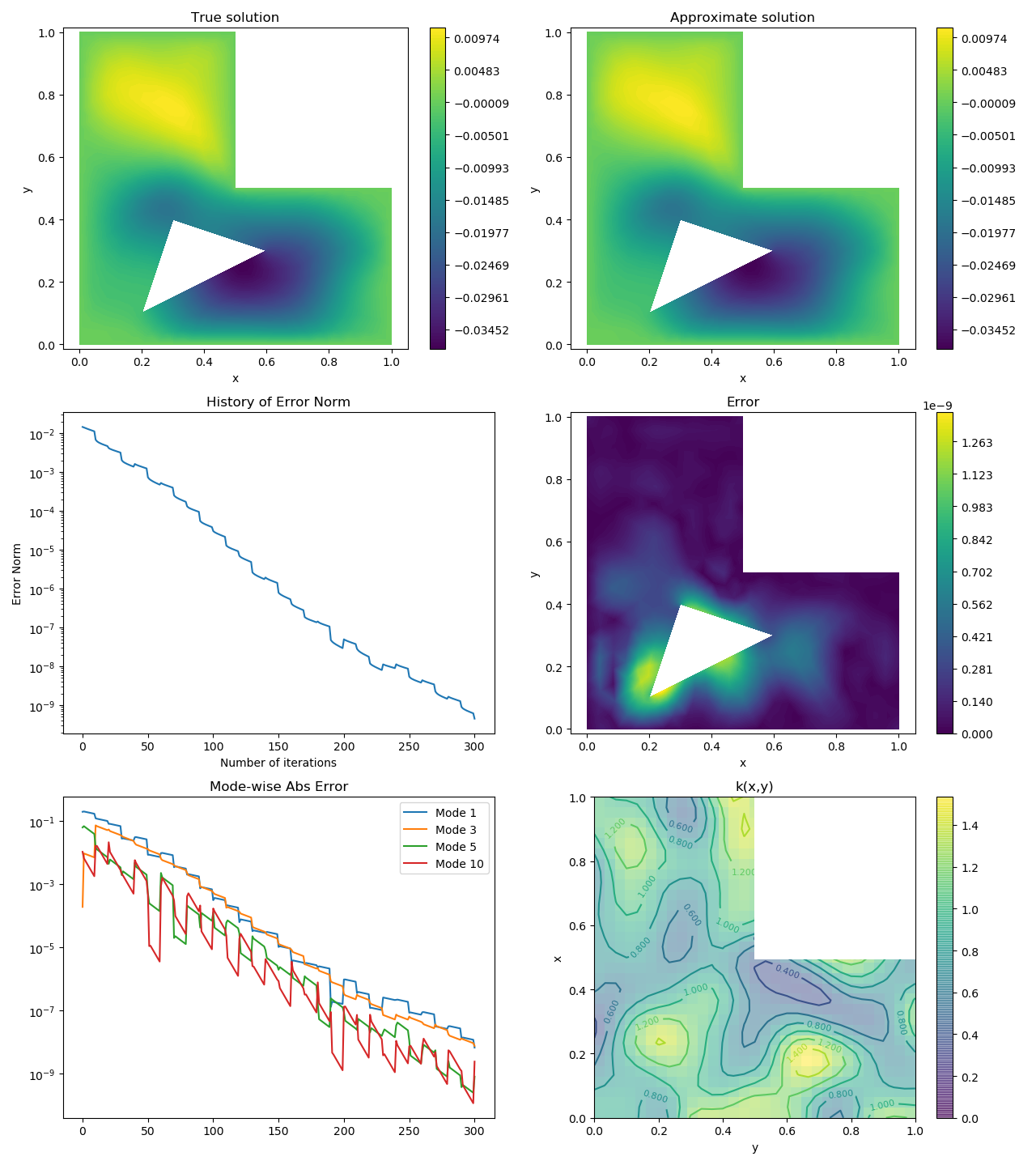}
         \caption{HINTS trained on $\Omega_L$, solving on $\Omega_L^{T2}$.}
         \label{fig:HINTS_triangle}
     \end{subfigure}
     \hfill
     \caption{Direct application of HINTS for domain adaptation. The results presented here are the plots of the target domain when HINTS has been directly used for inference without any re-training. Uniform convergence of the modes is observed as shown in the bottom left images for (a) and (b).}
     \label{fig:HINTS_geometries}
\end{figure}

\textbf{Direct Application:}\\ 
To begin with, we first investigate the domain adaptation capabilities of source DeepONet to make extrapolated approximations on target domain discretization for target HINTS. Precisely, the source DeepONet is directly employed and is iterated with the relaxation methods to approximate the solution of the dependent variable on the target discretization. On the target domains, we define $k(\boldsymbol{x})$, and $f(\boldsymbol x)$ by setting the input function values to be zero for points within the cutouts. The convergence of the solution is attributed to the generalization ability of DeepONet. The results for the two target tasks are presented in \autoref{fig:HINTS_geometries}. We observe that for Task1, convergence to machine precision is obtained in $182$ iterations. For Task2, even though much slower (takes longer than $300$ iterations), we do achieve convergence as well. We attribute the slow convergence to the three vertices of the triangle, which are singular points and are considered more difficult to handle. The HINTS was able to operate on this geometry and show good performance.

\textbf{Transfer Learning:}\\
As an alternative to the above approach, we propose to fine-tune the DeepONet with a small number of samples from the target domain. To that end, we employ the operator level transfer learning approach proposed in \cite{goswami2022deep}, where the target DeepONet is initialized with the weights and biases of the source DeepONet. To fine tune the target DeepONet, only the fully-connected layers of the branch network and the last layer of the trunk network are re-trained with a hybrid loss function that takes into account the difference in the conditional distribution of the source and the target. 

To train the target DeepONet for each of the two tasks, we generate $N_T = 500$ random samples of $k(\boldsymbol x)$ and $f(\boldsymbol x)$ and obtain the corresponding solution $u(\boldsymbol x)$. The target DeepONet is trained for $10\small{,}000$ iterations. When compared to the initial training of the network (from scratch), fine-tuning is faster. Once the target model is trained, we employ target DeepONet with the iterative solver on Target HINTS for the target domain.

\begin{figure}[!ht]
    \centering
    \includegraphics[width=12cm]{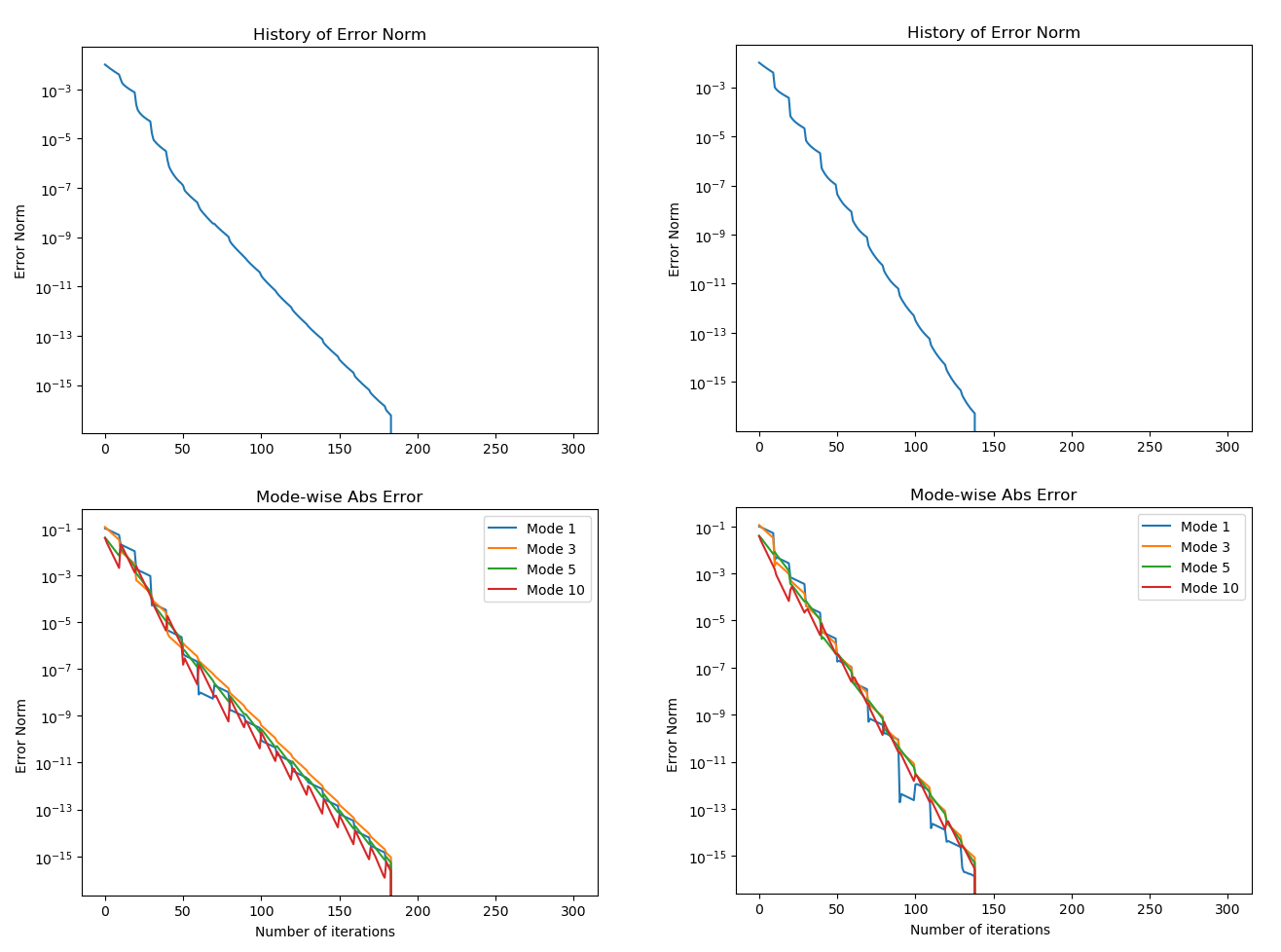}
    \caption{Convergence of the direct application of HINTS (left) and the transfer learning HINTS (right), done for the L-shaped Darcy problem with a circular cutout. The top figures show the error norm convergence and the bottom figures show the error norms of specific modes.}
    \label{fig:transfer_with_modes}
\end{figure}

\begin{figure}[!ht]
    \centering
    \includegraphics[width=10cm]{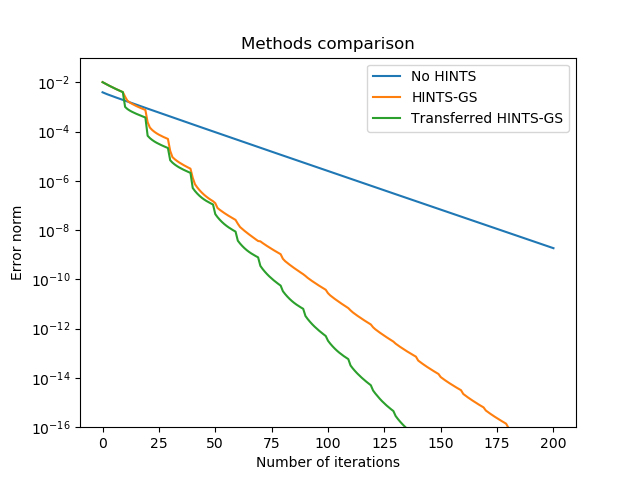}
    \caption{Convergence of the different methods for the L-shaped Darcy problem with a circular cutout. Comparison of the convergence of GS (without HINTS, shown in Blue Line), HINTS-GS (Orange Line) and the transfer learning HINTS-GS (Green Line) in terms of error decay.}
    \label{fig:compare_transfer}
\end{figure}

Finally, a comparative study is carried out based on the number of iterations each approach takes to converge to machine precision. We randomly select a sample from the target dataset and compare the convergence of the error (both the norm of the error and the mode errors) over iterations, between the direct application of HINTS and the transfer learning HINTS. This is shown in \autoref{fig:transfer_with_modes}. In addition, a comparison of the convergence rate that includes the standard GS solver is shown in \autoref{fig:compare_transfer}. We conclude that the HINTS solution is transferable to different domain geometries, and integrating the transfer-learning approach to fine tune the source DeepONet on the target domain results in a $25\%$ faster convergence without any loss of accuracy for this example. 

\subsection{Linear elasticity problem}

\begin{figure}[!ht]
     \centering
     \begin{subfigure}[t]{0.49\textwidth}
         \centering
         \includegraphics[width=\textwidth]{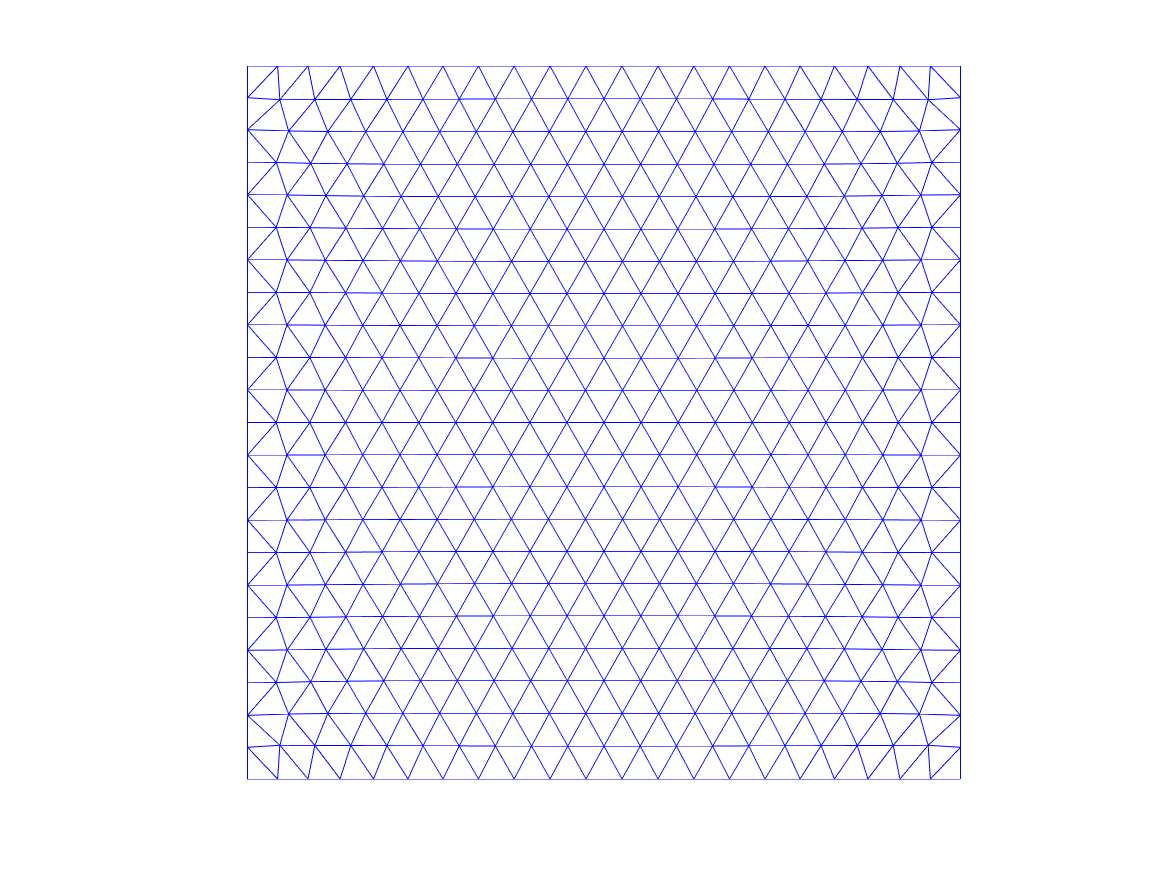}
         \caption{Square mesh ($\Omega^S$).}
     \end{subfigure}
     \hfill
     \begin{subfigure}[t]{0.49\textwidth}
         \centering
         \includegraphics[width=\textwidth]{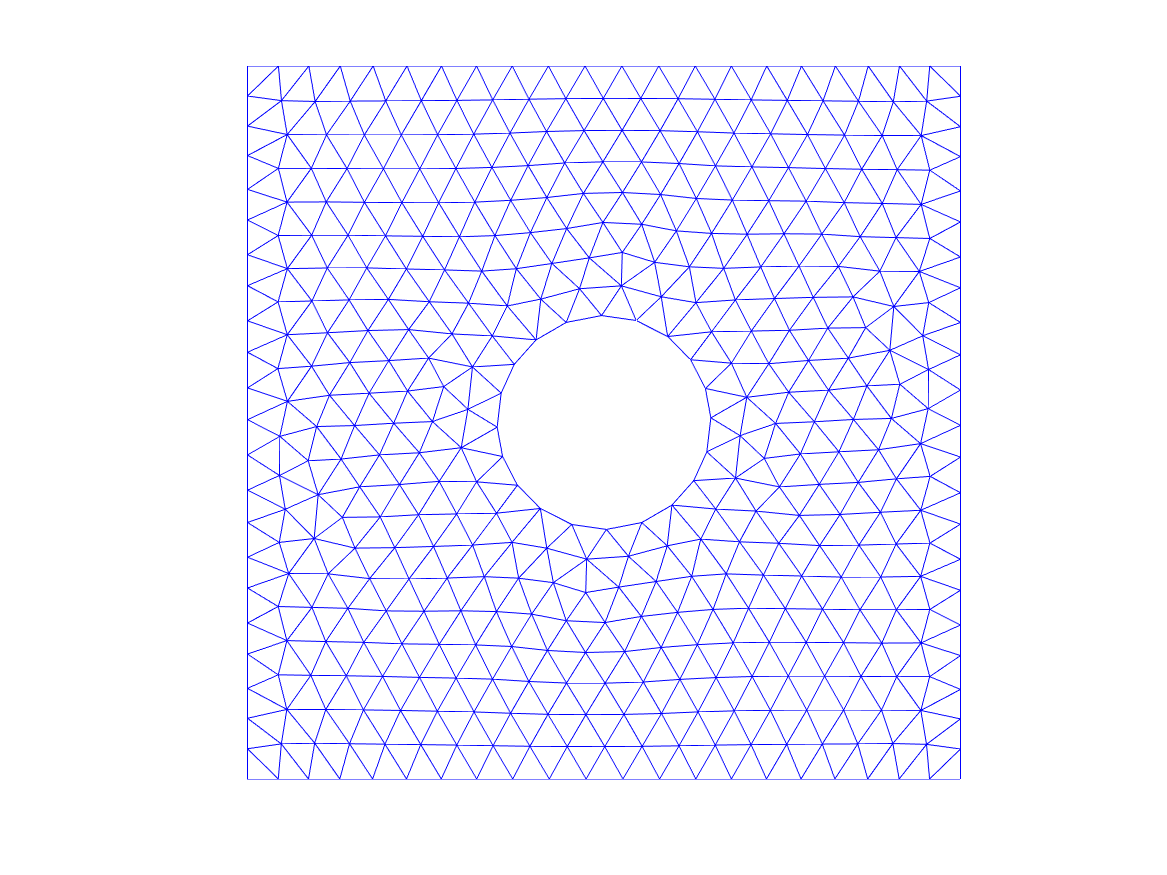}
         \caption{Square mesh with a circular hole ($\Omega^{T3}$).}
     \end{subfigure}
     \hfill
     \caption{Illustration of the meshes used for the elasticity numerical experiments.}
     \label{fig:square_domains}
\end{figure}

In the second example, we consider linear elasticity on a square domain under plane-strain conditions. The governing equation for the model is defined as:
\begin{align}
    \label{eqn:elasticity_PDE}
    \begin{cases}
        \epsilon_x(x, y) &= \frac{\partial u_x(x, y)}{\partial x} \\
        \epsilon_y(x, y) &= \frac{\partial u_x(x, y)}{\partial y} \\
        \gamma_{xy}(x, y) &= \frac{\partial u_x(x, y)}{\partial y} + \frac{\partial u_y(x, y)}{\partial x} \\
    \end{cases}
\end{align}
with subscripts $i,j\in\{1,2\}$ refer to the two in-plane directions, $\boldsymbol{u}_i$ is the displacement in the $i$ direction, $\varepsilon_{ij}$ is the strain component measured in $j$ direction due to displacement in $i$ direction, $\sigma_{ij}$ is the stress component, and $f_i$ is the body force in the $i$ direction. The Lamé parameters, $\mu=\frac{E}{2(1+\nu)}$ and $\lambda=\frac{\nu E}{(1+\nu)(1-2\nu)}$ describe the mechanical properties of the material, where $E$ and $\nu$ are the Young's modulus and the Poisson's ratio, respectively. In this example, we consider the square domain, $\Omega^S = [0, 1]\times [0, 1]$ as the source. The discretization of the domain is presented in \autoref{fig:square_domains}(a).

\bigbreak
\noindent
\textbf{Employing HINTS to solve the Elasticity problem on source domain}

The goal of the elasticity problem is to learn the operator of the system described by \autoref{eqn:elasticity_PDE}, which maps randomly generated spatially varying modulus of elasticity, $E(\boldsymbol x)$ and randomly varying force vector, $f(\boldsymbol{x})$ to the displacement vector, $\boldsymbol u (\boldsymbol x)$. In this example, we follow the same data generation approach as discussed in \autoref{subsec:darcy}, and the model was trained with $N_S = 85\small{,}000$ samples and tested with $N_S^{test} = 15\small{,}000$ samples. The branch net inputs two channels $(E(\boldsymbol{x})$ and $f(\boldsymbol x))$, and we adopt the same architecture for the convolution modules as discussed in the previous example. The fully-connected network following the convolution modules has a dimension $[256, 160]$. The dimension of the trunk network is $[2, 128, 128, 160]$. The network is trained for $25\small{,}000$ epochs to achieve roughly $4\%$  relative error, indicating a sufficiently well trained network. Now, we employ the HINTS algorithm to solve the elasticity problem on the source domain.

\textbf{Domain adaptation for the elasticity problem}

\begin{figure}[!ht]
     \centering
     \begin{subfigure}[t]{0.49\textwidth}
         \centering
         \includegraphics[width=\textwidth]{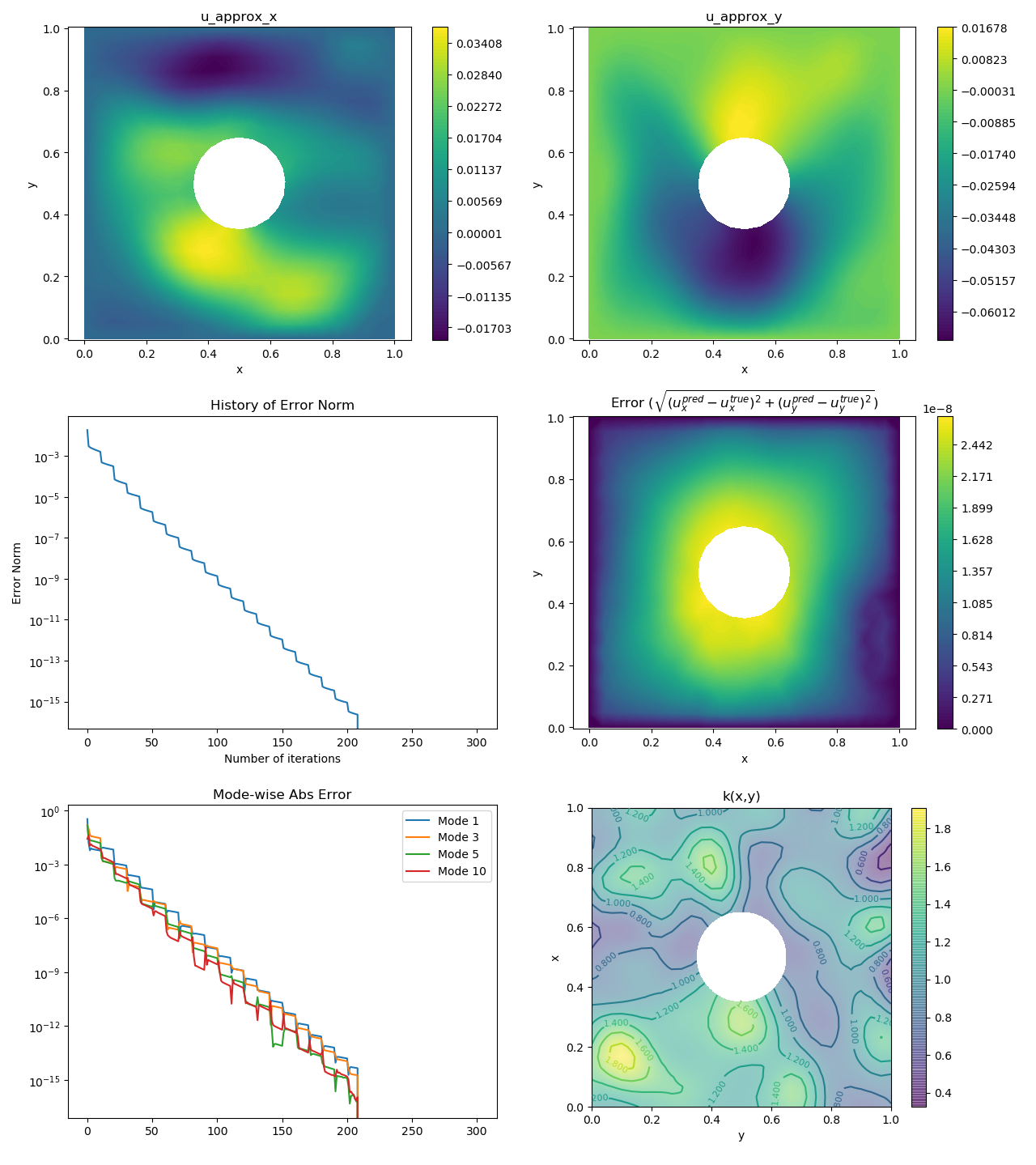}
         \caption{HINTS (direct application) trained on $\Omega^S$, solving on $\hat{\Omega}^S_{circle}$.}
         \label{fig:HINTS_elasticity_no_transfer}
     \end{subfigure}
     \hfill
     \begin{subfigure}[t]{0.49\textwidth}
         \centering
         \includegraphics[width=\textwidth]{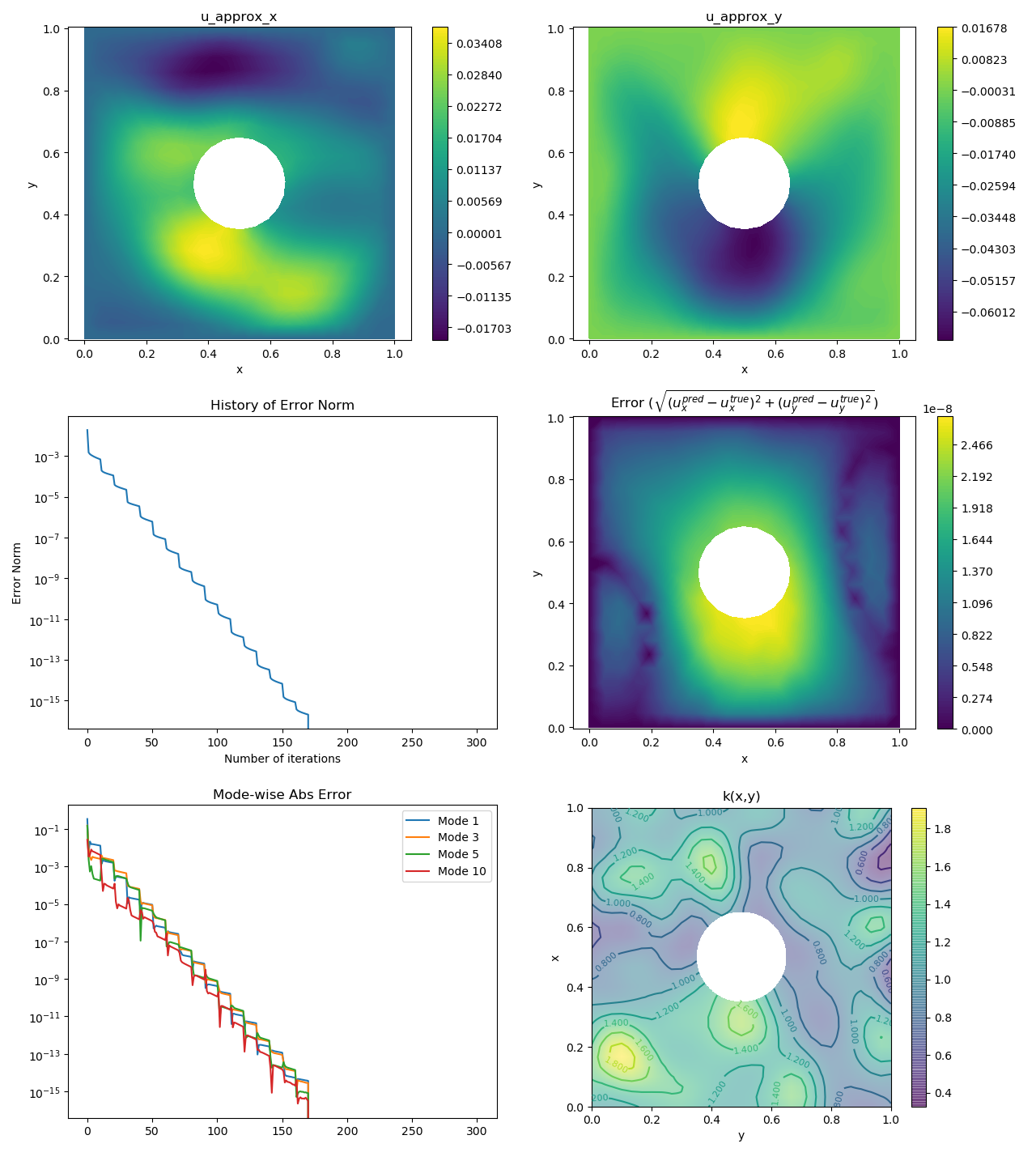}
         \caption{HINTS (transfer learning), solving on $\hat{\Omega}^S_{circle}$.}
         \label{fig:HINTS_elasticity_transfer}
     \end{subfigure}
     \hfill
     \caption{Results obtained using the direct application of HINTS (a) and transfer learning HINTS (b) for an example from the target data-set of the elasticity problem on the square domain with a circular cutout.}
     \label{fig:HINTS_elasticity}
\end{figure}

For investigating the domain adaptation capabilities, we define the following task:
\begin{itemize}
    \item \textbf{Task3:} From a square domain to a square domain with a circular cutout. The target domain is: $\Omega^{T3} = \Omega^S \backslash \{(x, y) | (x - 0.5) ^ 2 + (y - 0.5) ^ 2 \leq 0.15\}$.
\end{itemize}
The discretization of the target domains is shown in \autoref{fig:square_domains}. To begin with, we first employ the source DeepONet to infer on the target domain. In this scenario, no additional training is carried out on the target model. The HINTS algorithm is employed on the target domain for the iterative solver and use the source DeepONet to approximate the solution of displacement. The results presented in \autoref{fig:HINTS_elasticity_no_transfer}(a) show that using HINTS as is we can converge to machine precision after $205$ iterations on the target domain.  

As discussed in the previous example, we now integrate the operator level transfer learning algorithm with the HINTS model. In this setup, to train the target model, we generate samples of the input function by appending zeros to the function values with the circular cutout. The target model is trained with $N_T = 100$ samples, where the model is initialized with the optimized parameters of the source model and while training all the layers except the fully-connected layers of the branch network and the last layer of the trunk network are frozen. The fine tuned target DeepONet is replaced in the HINTS algorithm to generate the solution for the target domain. The results obtained using transfer learning integrate HINTS are presented in \autoref{fig:HINTS_elasticity_transfer}(b). In this setup, we observe convergence to machine precision after $169$ iterations, a $22.5\%$ improvement over the previous setup of employing HINTS with the source DeepONet to approximate solution for the target domain. 

\section{Summary}
\label{sec:conclusion}

In this work, we explored the transferability of HINTS on different computational domains for differential equations. Specifically, we considered two methods: (1) directly employ the HINTS for an equation defined in an unseen domain; (2) adopt the transfer operator learning to fine-tune the HINTS trained on the source domain for the unseen/target domain with limited data. Through presenting the results for Darcy flow and linear elasticity, we demonstrate the effectiveness of the two methods based on HINTS on fast, accurate solutions of differential equations. In particular, HINTS with transfer learning, by leveraging both the knowledge from the source domain and the target domain, converges even faster than using the direct application of HINTS on the target geometry. Despite the faster convergence, it comes with a price that one still needs a small dataset on the target domain. While we have focused on a specific instance of $k(\boldsymbol{x})$ ($E(\boldsymbol{x})$ for elasticity) and $f(\boldsymbol{x})$ in \autoref{sec:experiments}, the performance is consistent for the entire test dataset. \autoref{tab:quantified_results} shows the mean, median and standard deviation (STD) of the number of iterations needed for convergence for 100 different cases in the test dataset.

\begin{table}[htbp!]
    \centering
    \begin{tabular}{||c | c | c | c | c ||} 
        \hline
        Problem & Method & Mean & Median & STD \\ [0.5ex] 
        \hline\hline
        \multirow{3}{*}{Darcy} & GS & 403 & 405 & 34.62 \\ 
        \cline{2-5}
        & HINTS-GS & 165 & 164 & 19.03 \\
        \cline{2-5}
        & Transferred HINTS-GS & \textbf{157} & \textbf{162} & \textbf{15.75} \\
        \hline
        \multirow{3}{*}{Elasticity} & GS & 1029 & 1020 & 64.41 \\
        \cline{2-5}
        & HINTS-GS & 257 & 253 & 34.56 \\
        \cline{2-5}
        & Transferred HINTS-GS & \textbf{176} & \textbf{173} & \textbf{15.68} \\
        \hline
    \end{tabular}
    \caption{Summary of the results for the two benchmark problems. Mean, median and standard deviation of the number of iterations it takes for each method to converge to machine precision. The statistical measures were computed for 100 target samples of each one of the Darcy and elasticity problems.}
    \label{tab:quantified_results}
\end{table}

The capability of the direct application of HINTS for unseen geometries is, to some extent, rather surprising. Seemingly, a DeepONet trained for a fixed geometry (e.g., L-shaped domain) should not be effective on another geometry (e.g., L-shaped domain with a cutout) that is not included in the training dataset. We attribute the functionality of HINTS for unseen geometry to the following two factors: $(1)$ DeepONet simply needs to provide an approximate solution within HINTS, while the task of achieving accuracy is accomplished by the embedded numerical solver; $(2)$ the differential equation defined on the unseen geometry, for the examples that we consider, is similar to the equation defined on the original geometry but with an input function $k(\boldsymbol x)$ ($E(\boldsymbol x)$) defined in an extended domain. For $(1)$, intuitively, the prediction error of the DeepONet caused by the mismatch between the original and the new geometry depends on the difference between the two geometries. Within a reasonable degree of similarity between the two geometries, DeepONet can still decrease the errors of the low-frequency modes. For $(2)$, using the case of Darcy flow as an example, it may be shown that the differential equation defined in the L-shaped domain excluding the cutout is equivalent to the same equation defined in the L-shaped domain, where (a) within the cutout $k(\boldsymbol x)$ is simply padded with zero, and (b) the boundary condition at the cutout boundary is zero Neumann boundary condition. Technically, our approach of padding $k(\boldsymbol x)$ inside the cutout with zeros conforms with such equivalence. Therefore, generalizing the L-shaped domain into a new geometry (L-shaped domain with a cutout) is transformed into the generalization of $k(\boldsymbol x)$ from GRF in the training dataset into an unseen $k(\boldsymbol x)$, where it is from GRF outside the cutout but equals to zero inside the cutout.

\bibliographystyle{unsrt}  
\bibliography{references}

\begin{thebibliography}{10}

\bibitem{hughes2012finite}
Thomas~JR Hughes.
\newblock {\em The finite element method: linear static and dynamic finite
  element analysis}.
\newblock Courier Corporation, 2012.

\bibitem{simo2006computational}
Juan~C Simo and Thomas~JR Hughes.
\newblock {\em Computational inelasticity}, volume~7.
\newblock Springer Science \& Business Media, 2006.

\bibitem{hughes2005isogeometric}
Thomas~JR Hughes, John~A Cottrell, and Yuri Bazilevs.
\newblock Isogeometric analysis: {CAD}, finite elements, {NURBS}, exact
  geometry and mesh refinement.
\newblock {\em Computer Methods in Applied Mechanics and Engineering},
  194(39-41):4135--4195, 2005.

\bibitem{jing2002numerical}
Lanru Jing and JA~Hudson.
\newblock Numerical methods in rock mechanics.
\newblock {\em International Journal of Rock Mechanics and Mining Sciences},
  39(4):409--427, 2002.

\bibitem{rappaz2003numerical}
Michel Rappaz, Michel Bellet, Michel~O Deville, and R~Snyder.
\newblock {\em Numerical modeling in materials science and engineering}.
\newblock Springer, 2003.

\bibitem{goswami2019adaptive}
Somdatta Goswami, Cosmin Anitescu, and Timon Rabczuk.
\newblock Adaptive phase field analysis with dual hierarchical meshes for
  brittle fracture.
\newblock {\em Engineering Fracture Mechanics}, 218:106608, 2019.

\bibitem{bharali2022robust}
Ritukesh Bharali, Somdatta Goswami, Cosmin Anitescu, and Timon Rabczuk.
\newblock A robust monolithic solver for phase-field fracture integrated with
  fracture energy based arc-length method and under-relaxation.
\newblock {\em Computer Methods in Applied Mechanics and Engineering},
  394:114927, 2022.

\bibitem{zhang2022g2}
Enrui Zhang, Bart Spronck, Jay~D Humphrey, and George~Em Karniadakis.
\newblock G2{$\Phi$}net: Relating genotype and biomechanical phenotype of
  tissues with deep learning.
\newblock {\em arXiv preprint arXiv:2208.09889}, 2022.

\bibitem{goswami2022neural}
Somdatta Goswami, David~S Li, Bruno~V Rego, Marcos Latorre, Jay~D Humphrey, and
  George~Em Karniadakis.
\newblock Neural operator learning of heterogeneous mechanobiological insults
  contributing to aortic aneurysms.
\newblock {\em arXiv preprint arXiv:2205.03780}, 2022.

\bibitem{patera1984spectral}
Anthony~T Patera.
\newblock A spectral element method for fluid dynamics: laminar flow in a
  channel expansion.
\newblock {\em Journal of Computational Physics}, 54(3):468--488, 1984.

\bibitem{kim1987turbulence}
John Kim, Parviz Moin, and Robert Moser.
\newblock Turbulence statistics in fully developed channel flow at low
  {R}eynolds number.
\newblock {\em Journal of Fluid Mechanics}, 177:133--166, 1987.

\bibitem{cockburn2012discontinuous}
Bernardo Cockburn, George~E Karniadakis, and Chi-Wang Shu.
\newblock {\em Discontinuous {G}alerkin methods: theory, computation and
  applications}, volume~11.
\newblock Springer Science \& Business Media, 2012.

\bibitem{zhang2022hybrid}
Enrui Zhang, Adar Kahana, Eli Turkel, Rishikesh Ranade, Jay Pathak, and
  George~Em Karniadakis.
\newblock A hybrid iterative numerical transferable solver (hints) for pdes
  based on deep operator network and relaxation methods.
\newblock {\em arXiv preprint arXiv:2208.13273}, 2022.

\bibitem{goswami2022deep}
Somdatta Goswami, Katiana Kontolati, Michael~D Shields, and George~Em
  Karniadakis.
\newblock Deep transfer learning for partial differential equations under
  conditional shift with deeponet.
\newblock {\em arXiv preprint arXiv:2204.09810}, 2022.

\bibitem{briggs2000multigrid}
William~L Briggs, Van~Emden Henson, and Steve~F McCormick.
\newblock {\em A multigrid tutorial}.
\newblock SIAM, 2000.

\bibitem{bramble2019multigrid}
James~H Bramble.
\newblock {\em Multigrid methods}.
\newblock Chapman and Hall/CRC, 2019.

\bibitem{van2007spectral}
Martin~B van Gijzen, Yogi~A Erlangga, and Cornelis Vuik.
\newblock Spectral analysis of the discrete helmholtz operator preconditioned
  with a shifted laplacian.
\newblock {\em SIAM Journal on Scientific Computing}, 29(5):1942--1958, 2007.

\bibitem{raissi2019physics}
Maziar Raissi, Paris Perdikaris, and George~E Karniadakis.
\newblock Physics-informed neural networks: A deep learning framework for
  solving forward and inverse problems involving nonlinear partial differential
  equations.
\newblock {\em Journal of Computational physics}, 378:686--707, 2019.

\bibitem{goswami2020transfer}
Somdatta Goswami, Cosmin Anitescu, Souvik Chakraborty, and Timon Rabczuk.
\newblock Transfer learning enhanced physics informed neural network for
  phase-field modeling of fracture.
\newblock {\em Theoretical and Applied Fracture Mechanics}, 106:102447, 2020.

\bibitem{moore2022learning}
Nicholas~S Moore, Eric Cyr, and Christopher Siefert.
\newblock Learning an algebriac multrigrid interpolation operator using a
  modified graphnet architecture.
\newblock Technical report, Sandia National Lab.(SNL-NM), Albuquerque, NM
  (United States), 2022.

\bibitem{luz2020learning}
Ilay Luz, Meirav Galun, Haggai Maron, Ronen Basri, and Irad Yavneh.
\newblock Learning algebraic multigrid using graph neural networks.
\newblock In {\em International Conference on Machine Learning}, pages
  6489--6499. PMLR, 2020.

\bibitem{gotz2018machine}
Markus G{\"o}tz and Hartwig Anzt.
\newblock Machine learning-aided numerical linear algebra: convolutional neural
  networks for the efficient preconditioner generation.
\newblock In {\em 2018 IEEE/ACM 9th Workshop on Latest Advances in Scalable
  Algorithms for Large-Scale Systems (scalA)}, pages 49--56. IEEE, 2018.

\bibitem{goswami2022operator}
Somdatta Goswami, Aniruddha Bora, Yue Yu, and George~Em Karniadakis.
\newblock Physics-informed neural operators.
\newblock {\em arXiv preprint arXiv:2207.05748}, 2022.

\bibitem{li2021fourier}
Zongyi Li, Nikola~Borislavov Kovachki, Kamyar Azizzadenesheli, Burigede liu,
  Kaushik Bhattacharya, Andrew Stuart, and Anima Anandkumar.
\newblock {Fourier Neural Operator for Parametric Partial Differential
  Equations}.
\newblock In {\em In Proceedings of the International Conference on Learning
  Representations}, 2021.

\bibitem{lu2021learning}
Lu~Lu, Pengzhan Jin, Guofei Pang, Zhongqiang Zhang, and George~Em Karniadakis.
\newblock Learning nonlinear operators via deeponet based on the universal
  approximation theorem of operators.
\newblock {\em Nature Machine Intelligence}, 3(3):218--229, 2021.

\bibitem{goswami2022physics}
Somdatta Goswami, Minglang Yin, Yue Yu, and George~Em Karniadakis.
\newblock A physics-informed variational deeponet for predicting crack path in
  quasi-brittle materials.
\newblock {\em Computer Methods in Applied Mechanics and Engineering},
  391:114587, 2022.

\bibitem{kontolati2022influence}
Katiana Kontolati, Somdatta Goswami, Michael~D Shields, and George~Em
  Karniadakis.
\newblock On the influence of over-parameterization in manifold based
  surrogates and deep neural operators.
\newblock {\em arXiv preprint arXiv:2203.05071}, 2022.

\bibitem{kingma2014adam}
Diederik~P Kingma and Jimmy Ba.
\newblock {Adam: A method for stochastic optimization}.
\newblock {\em arXiv preprint arXiv:1412.6980}, 2014.

\bibitem{paszke2019pytorch}
Adam Paszke, Sam Gross, Francisco Massa, Adam Lerer, James Bradbury, Gregory
  Chanan, Trevor Killeen, Zeming Lin, Natalia Gimelshein, Luca Antiga, et~al.
\newblock Pytorch: {A}n imperative style, high-performance deep learning
  library.
\newblock {\em Advances in neural information processing systems},
  32:8026--8037, 2019.

\end{thebibliography}

\end{document}